\documentclass[draft,notitlepage]{article}
 
\usepackage{amsmath}
\usepackage{amsthm}
 \usepackage[Symbol]{upgreek}
 \usepackage{amssymb}

\DeclareMathAlphabet{\mathbf}{OML}{cmm}{b}{it}

\numberwithin{equation}{section}

\newcommand{\bbold}{\mathbb}
\renewcommand\epsilon{\varepsilon}

\def\R { {\bbold R} }
\def\Q { {\bbold Q} }

\def\N { {\bbold N} }

\def\da {\operatorname{da}}

\def\divides {\,|\,}
\newcommand{\udots}{{\mathinner{\mskip1mu\raise1pt\vbox{\kern7pt\hbox{.}}\mskip2mu\raise4pt\hbox{.}\mskip2mu\raise7pt\hbox{.}\mskip1mu}}}
\newcommand{\nobracket}{}
\newcommand{\nosymbol}{}

\def \ex{\operatorname{e}}

\def \Cc{\mathcal{C}}

\def \d{\operatorname{d}}

\def \<{\langle}
\def \>{\rangle}

\def \tilde {\widetilde}

\def \hat {\widehat}

\def \supp {\operatorname{supp}}

\def \((  {(\!(}
\def \)) {)\!)}
\def \T{\mathbb{T}}

\DeclareMathSymbol{\precequ}{\mathrel}{symbols}{"16}
\DeclareMathSymbol{\succequ}{\mathrel}{symbols}{"17}

\def \nasymp{\not\asymp}

\newtheorem{theorem}{Theorem}[section]

\newtheorem{cor}[theorem]{Corollary}

\theoremstyle{definition}

\newtheorem{definition}[theorem]{Definition}

\theoremstyle{remark}
\newtheorem{conjecture}[theorem]{Conjecture}

\newcommand{\abs}[1]{\lvert#1\rvert}

\def \fM {{\mathfrak M}}

\def \fm {{\mathfrak m}}
\def \fn {{\mathfrak n}}

\def \fM {{\mathfrak M}}

\def \fm {{\mathfrak m}}
\def \fn {{\mathfrak n}}

\let\oldi\i
\let\oldj\j
\renewcommand\i{\relax\ifmmode{\boldsymbol{i}}\else\oldi\fi}
\renewcommand\j{\relax\ifmmode{\boldsymbol{j}}\else\oldj\fi}

\renewcommand\leq{\leqslant}
\renewcommand\geq{\geqslant}
\renewcommand\preceq{\preccurlyeq}
\renewcommand\succeq{\succcurlyeq}

\DeclareFontFamily{U}{fsy}{}
\DeclareFontShape{U}{fsy}{m}{n}{<->s*[.9]psyr}{}
\DeclareSymbolFont{der@m}{U}{fsy}{m}{n}
\DeclareMathSymbol{\der}{\mathord}{der@m}{182}

\DeclareSymbolFont{der@m}{U}{fsy}{m}{n}
\DeclareMathSymbol{\derdelta}{\mathord}{der@m}{100}

\makeatletter
\newcommand{\raisemath}[1]{\mathpalette{\raisem@th{#1}}}
\newcommand{\raisem@th}[3]{\raisebox{#1}{$#2#3$}}
\makeatother

\DeclareSymbolFont{imag@m}{OT1}{cmr}{m}{ui}
\DeclareMathSymbol{\imag}{\mathord}{imag@m}{105}

\DeclareFontFamily{OMS}{smallo}{}
\DeclareFontShape{OMS}{smallo}{m}{n}{<->s*[.65]cmsy10}{}
\DeclareSymbolFont{smallo@m}{OMS}{smallo}{m}{n}
\DeclareMathSymbol{\smallo}{\mathord}{smallo@m}{79}

\DeclareFontFamily{OMS}{largerdot}{}
\DeclareFontShape{OMS}{largerdot}{m}{n}{<->s*[.8]cmsy10}{}
\DeclareSymbolFont{largerdot@m}{OMS}{largerdot}{m}{n}
\DeclareMathSymbol{\largerdot}{\mathord}{largerdot@m}{15}

\DeclareMathSymbol{\llambda}{\mathord}{der@m}{108}
\DeclareMathSymbol{\rrho}{\mathord}{der@m}{114}

\def \upl{\uplambda}
\def \Upl{\Uplambda}
\def \upo{\upomega}
\def \Upo{\Upomega}

\newcommand{\equationqed}[1]{\[\pushQED{\qed}#1 \qedhere\popQED\]\let\qed\relax}
\newcommand{\alignqed}[1]{\begin{align*}\pushQED{\qed} #1 \qedhere\popQED\end{align*}\let\qed\relax}

\makeatletter
\everydisplay=\expandafter{\the\everydisplay
  \ifdim\@totalleftmargin>\z@
    \displayindent=\dimexpr+\leftmargin\relax
    \displaywidth=\columnwidth
  \fi
}
\makeatother

\renewenvironment{abstract}
 {\small
  \begin{center}
  \bfseries \abstractname\vspace{-.5em}\vspace{0pt}
  \end{center}
  \list{}{
    \setlength{\leftmargin}{2.5em}%
    \setlength{\rightmargin}{\leftmargin}%
  }%
  \item\relax}
 {\endlist}

\begin{document}

\title{On Numbers, Germs, and Transseries}

%\author{Matthias Aschenbrenner}
%\address{Department of Mathematics\\
%University of California, Los Angeles\\
%Los Angeles, CA 90095\\
%U.S.A.}
%\email{matthias@math.ucla.edu}

%\author{Lou van den Dries}
%\address{Department of Mathematics\\
%University of Illinois at Urbana-Cham\-paign\\
%Urbana, IL 61801\\
%U.S.A.}
%\email{vddries@math.uiuc.edu}

%\author{Joris van der Hoeven}
%\address{\'Ecole Polytechnique\\
%91128 Palaiseau Cedex\\
%France}
%\email{vdhoeven@lix.polytechnique.fr}

\author{\normalsize Matthias Aschenbrenner, Lou van den Dries, Joris van der Hoeven}

\date{}

%Redefine footnote temporarily to produce MSC
\renewcommand{\thefootnote}{\fnsymbol{footnote}} 
\footnotetext{\emph{Date}. November 16, 2017}     
\footnotetext{The first-named author was partially support by NSF Grant DMS-1700439.}     
\renewcommand{\thefootnote}{\arabic{footnote}}

\maketitle

\begin{abstract}
\noindent
Germs of real-valued functions, surreal numbers, and transseries
are three ways to enrich the real continuum by infinitesimal and infinite quantities. 
%surreal numbers have a combinatorial flavor and encompass all ordinals;  germs of (differentiable) real-valued functions are central objects in asymptotic analysis;  and transseries are formal representations for the growth rates of such germs at infinity.
Each of these comes with naturally interacting notions of {\it ordering}\/ and   {\it derivative}.
The category of $H$-fields provides a common framework for the relevant  
algebraic structures.
We give an exposition of our results on the model theory of $H$-fields, and we report on recent progress in unifying germs, surreal numbers, and transseries from the point of view of asymptotic differential\nolinebreak{}~algebra.
\end{abstract}

\medskip

\noindent
%It was a fundamental discovery of G.~Cantor that his theory of sets, for the first time, allowed the rigorous treatment of ``actual infinite'' quantities, as opposed to the ``potential infinity'' exhibited, say, by sequences of real numbers. Another success of his theory was to give   mathematically sound constructions of the real numbers, resulting in a continuum from which infinitesimals are banished. One such construction, by R.~Dedekind, introduces the reals using his eponymous cuts of rationals.
Contemporaneous with Can\-tor's work in the 1870s but less well-known, P.~du~Bois-Reymond~\cite{dBR71}--\cite{dBR82} had original ideas concerning  non-Cantorian infinitely large and small quantities~\cite{Ehr06}. 
%an early pioneer of this subject, %Motivated by investigations of convergence criteria for infinite series,
He developed a ``calculus of infinities'' to deal with the growth rates of functions of one real variable, %approaching $+\infty$, 
representing their ``potential infinity'' by an ``actual infinite'' quantity.  The reciprocal of a function tending to infinity is one which tends to zero, hence represents an ``actual infinitesimal''. 

These ideas were unwelcome to Cantor~\cite{Fisher} and misunderstood by him,
but were made rigorous by F.~Hausdorff~\cite{Hau06}--\cite{Hau09} and G.~H.~Hardy~\cite{Har10}--\cite{Har13}. 
Hausdorff firmly grounded du~Bois-Reymond's ``orders of infinity'' in Can\-tor's   set-theoretic universe~\cite{Felgner}, while
Hardy focused on their differential aspects and introduced the
{\it logarithmico-exponential functions}\/
(short: {\it LE-functions}). This led  to the  concept of a {\it Hardy field}\/ (Bour\-ba\-ki~\cite{Bour51}), developed further mainly by Rosenlicht~\cite{Ros83}--\cite{Ros95} and Boshernitzan~\cite{Bos81}--\cite{Bos87}.   For the role of Hardy fields in {\em o-minimality} see~\cite{Miller}.

{\it Surreal numbers}\/ were discovered (or created?) in the 1970s by J.~H.~Conway~\cite{Con76} 
%while studying combinatorial games, and
%(The title of our paper is a play on that of Conway's classic~\cite{Con76}.) 
and popularized by M.~Gardner, and by D.~E.~Knuth~\cite{Knuth} who coined the term ``surreal number''. The surreal numbers form a proper class containing all reals as well as Cantor's ordinals, and come equipped with a natural ordering and arithmetic operations turning them into an ordered field. 
Thus with $\omega$ the first infinite ordinal, $\omega-\pi$, $1/\omega$, $\sqrt{\omega}$  make sense as surreal numbers. In contrast to non-standard real numbers, their construction is completely canonical, naturally generalizing both Dedekind cuts and von Neumann's construction of the ordinals.
 (In the words of  their creator~\cite[p.~102]{Con94}, the surreals 
are %``are nothing less than 
``the only correct extension of the notion of real number to the infinitely large and the infinitesimally small.'') The surreal universe is very rich,   yet 
shares many properties with the real world.
For example, the ordered field of surreals is real closed  and hence, by Tarski~\cite{Tarski}, 
an elementary extension of its ordered subfield of real numbers.
(In fact, every set-sized real closed field embeds into the field of surreal numbers.)
M.~Kruskal anticipated the  use of surreal numbers in asymptotics, and based on his ideas Gonshor~\cite{Gon86}  extended the  exponential function on the reals to
one on the surreals, with the same first-order logical properties~\cite{DE}. 
Rudiments of analysis for functions on the surreal numbers
have also been developed~\cite{All87,CEF15,RSS14}.

{\it Transseries}\/ generalize LE-functions in a similar way that surreals generalize reals and ordinals. Transseries have a  precursor in the {\em generalized power series\/} of Levi-Civita~\cite{LC92, LC98}  and Hahn~\cite{Hahn}, but were only systematically considered in the 1980s, independently by \'Ecalle~\cite{Ec92} and Dahn-G\"oring~\cite{DG}. 
\'Ecalle introduced transseries as formal counterparts to his ``analyzable functions'',  which were central to
his work on Dulac's Problem (related to  Hil\-bert's~16th Problem on polynomial vector fields). 
Dahn and G\"oring were motivated by Tarski's Problem on the model theory of the real field with
exponentiation.
Transseries have since been used in various parts of mathematics  and physics; their formal nature also makes them suitable for calculations in computer algebra systems.
Key examples of transseries are the {\it logarithmic-exponential series}\/ 
({\it LE-series}\/ for short)~\cite{DMM97,DMM01}; more general notions of transseries have
been introduced in~\cite{vdH:phd,Schm01}.
%The {\it raison d'\^{e}tre}\/ of 
A transseries can represent a function of a real variable using exponential and logarithmic terms, going beyond the more prevalent 
asymptotic expansions in terms of powers of the independent variable. Transseries can  be manipulated algebraically---added, subtracted, multiplied, divided---and like power series,  can be differentiated term-wise: they comprise a differential field.
However, they carry much more structure: for example, by virtue of its construction,~the field of LE-series comes with an  exponential function;
there is a natural notion of composition for transseries; 
and differential-compositional equations in transseries are sometimes amenable to functional-analytic techniques~\cite{vdH:noeth}.

The logical properties of the {\it exponential\/} field of LE-series have been well-un\-der\-stood since the 1990s:
by \cite{Wilkie} and \cite{DMM97} it is model-complete and o-minimal.
In our book~\cite{vdH:b} we focused instead on the {\it differential}\/  field of LE-series, denoted below by $\T$,
and obtained some decisive results about its model theory. 
Following A.~Robinson's general ideas we 
placed $\T$ into a suitable category of \textit{$H$-fields}\/ and, by developing the extension theory of $H$-fields, showed  that $\T$ is existentially closed as an
$H$-field: each system of algebraic differential equations and inequalities over~$\T$ which has a solution in an $H$-field extension of~$\T$ already has one in $\T$ itself. 
In~\cite{vdH:b} we also prove the related fact that $\T$ is model-complete; indeed, we obtain a quantifier elimination (in a natural language) for $\mathbb T$.
As a consequence, the elementary theory of $\mathbb T$ is decidable, and  model-theoretically ``tame''
in various ways: for example, it has Shelah's {\it non-independence property} (NIP). 

Results from~\cite{vdH:b} about existential closedness, model completeness, and quantifier
elimination substantiate the
intuition, expressed already in~\cite{Ec92}, that~$\mathbb T$  plays the role of a {\em universal domain}\/ for the part of
asymptotic differential algebra that steers clear of oscillations.
How far does this intuition lead us? Hardy's field of LE-functions embeds into $\mathbb{\T}$, as an ordered differential field, but 
this fails for other Hardy fields. The natural question here is:  {\it Are all maximal Hardy fields
elementarily equivalent to $\mathbb{T}$?}\/ 
It would mean that any maximal Hardy field instantiates Hardy's vision of a maximally inclusive and well-behaved algebra of oscillation-free real functions. Related is
the issue of embedding Hardy fields into more general differential fields of transseries. 
Positive answers to these questions  would
tighten the link between germs of functions (living in Hardy fields) and their transseries expansions. 
We may also ask how surreal numbers fit into the picture: {\it Is there a natural isomorphism between the field of surreal numbers and some field  of generalized transseries?}\/ This would   make it possible to differentiate and compose surreal numbers as if they were functions, and   confirm Kruskal's premonition of a connection between surreals and the asymptotics of functions.

We believe that answers to these questions are within grasp due to advances in our understanding 
%of the ideas of number, series, and function during the last decade, and would exhibit a heretofore unknown relationship between them. 
during the last decade as represented in our book~\cite{vdH:b}.
We discuss these questions with more details in Sections~\ref{HH},~\ref{HT},~\ref{HN}.
In Section~\ref{sec:trans} we set the stage by describing  Hardy fields and transseries as two competing approaches to the asymptotic behavior of non-oscillatory real-valued functions. 
(Section~\ref{HN} includes a brief synopsis of the remarkable  surreal number system.) 
In Section~\ref{sec:H-fields} we define  $H$-fields and state the main results of~\cite{vdH:b}.

%Can transseries serve as a bridge between germs and numbers? 

We let $m$, $n$ range over $\N=\{0,1,2,\dots\}$.
Given an (additive) abelian group $A$ we let $A^{\neq}:=A\setminus\{0\}$. In some places below we assume familiarity with very basic model theory, for example,
on the level of \cite[Appendix~B]{vdH:b}. ``Definable'' will mean ``definable with parameters''.

\section{Orders of Infinity and Transseries}\label{sec:trans}

\subsection*{Germs of continuous functions}

\noindent
Consider continuous real-valued functions  whose domain is a subset of $\R$ containing an interval~$(a,+\infty)$,  $a\in\R$. Two such  functions {\it have the same germ}\/ (at~$+\infty$) if
they agree on an interval $(a,+\infty)$,  $a\in\R$, contained in both their domains; this defines an equivalence relation on the set of
such functions, whose equivalence classes are called {\it germs.}\/ 
Addition and multiplication of germs is defined pointwise, giving rise to a commutative ring~$\mathcal C$.
 For a germ $g$ of such a function we also let $g$ denote that function if the resulting ambiguity is harmless.
With this convention, given a property~$P$ of real numbers
and~$g\in\mathcal C$ we say that~$P\big(g(t)\big)$ holds {\it eventually}\/ if~$P\big(g(t)\big)$ holds for all sufficiently large real $t$ in the domain of~$g$. 
%For a property~$F$ of functions $(a,+\infty)\to\R$, $a\in\R$, which is preserved under restriction
%to each subinterval $(b,+\infty)$ ($b\in\R$, $b\geq a$), we also say that a germ has property $F$
%if it has a representative with this property.
%We let $\mathcal{C}$ be the subring of $\mathcal{G}$ consisting of continuous germs. 
 We identify each real number $r$ with the germ  of the constant function 
$\R\to \R$ with value $r$. This  makes the field $\R$ into a subring of $\mathcal{C}$.

Following Hardy we define for $f,g\in\mathcal{C}$, 
\begin{align*} f\preceq g \quad:&\Longleftrightarrow\quad \text{for some $c\in\R^{>0}$  we have
$\abs{f(t)}\leq c\abs{g(t)}$ eventually},\\
f \prec g \quad:&\Longleftrightarrow\quad \text{for every $c\in\R^{>0}$  we have
$\abs{f(t)} <  c\abs{g(t)}$ eventually}.
\end{align*} 
The reflexive and transitive relation $\preceq$ yields an equivalence relation
$\asymp$ on~$\mathcal{C}$  by setting
$f\asymp g :\Longleftrightarrow f\preceq g$ and $g\preceq f$,
%$$f\asymp g \quad:\Longleftrightarrow\quad \text{$f\preceq g$ and $g\preceq f$,}$$
and $\preceq$ induces a partial ordering on the set of equivalence classes of~$\asymp$; these equivalence classes are
essentially du~Bois-Reymond's ``orders of infinity''.  
%Also, $f\prec g\Longrightarrow f\preceq g \text{ and }f\not\asymp g$.  
Thus with  $x$ the germ of the identity function\nopagebreak{}~on~$\R$:
$$0\ \prec\ 1\ \prec\ \log\log x\ \prec\ \log x\ \prec\ \sqrt{x}\ \prec\  x\ \asymp\ -2x + x\sin x\ \prec\ x^2\ \prec\ \ex^x.$$ 
One way to create interesting subrings of $\mathcal C$ is via 
  expansions of the  field of real numbers: any such  expansion~$\tilde{\R}$
  %$\mathbf R$ 
 gives rise to the subring
%$H(\mathbf R)$ 
$H(\tilde{\R})$ of $\mathcal C$ consisting of
the germs of the continuous functions~$\R\to\R$ that are definable in $\tilde{\R}$.
%$\mathbf R$. 

\subsection*{Hausdorff fields}

\noindent
A {\em Hausdorff field\/} is by definition a subfield of $\mathcal{C}$. 
Simple examples  are 
\begin{equation}\label{ex:Hausdorff fields}
\Q,\quad \R,\quad
\R(x), \quad \R(\sqrt{x}),\quad \R(x,\ex^x,\log x).
\end{equation}
That $\R(x,\ex^x,\log x)$ is a Hausdorff field,  for instance, 
follows from two easy facts: first, an element $f$ of $\mathcal{C}$ is a unit
iff $f(t)\neq 0$ eventually (and then either~${f(t)>0}$ eventually or $f(t)<0$ eventually), and if
$f\neq 0$ is an element of the sub\-ring $\R[x,\ex^x,\log x]$ of~$\mathcal{C}$, then $f\asymp x^k \ex^{l x}(\log x)^m$
for some  $k,l,m\in\N$.
Alternatively, one can use the fact that  an expansion $\tilde{\R}$
%$\mathbf R$ 
of the field of reals is o-minimal iff~$H(\tilde{\R})$ is a Hausdorff field, and note that
the examples above are subfields of~$H(\R_{\exp})$ where~$\R_{\exp}$ is the exponential field of real numbers, which is well-known to be o-minimal by Wilkie~\cite{Wilkie}. 

Let $H$ be a Hausdorff field. Then $H$ becomes an ordered field with (total) ordering given by:
$f >0$ iff $f(t) >0$ eventually.
Moreover, the set of orders of infinity in~$H$ is totally ordered by $\preceq$: for $f,g\in H$
we have $f\preceq g$ or $g\preceq f$. In his landmark paper~\cite{Hau09},
Hausdorff  essentially proved that $H$ has a unique algebraic Hausdorff field extension that is
real closed. (Writing before Artin and Schreier~\cite{AS},  of course he doesn't use this terminology.)
He was particularly interested in ``maximal'' objects and their order type. 
By Hausdorff's Maximality Principle (a form of Zorn's Lemma) every Hausdorff field is contained in one that is maximal with respect to inclusion.
By the above, maximal Hausdorff fields are real closed.
Hausdorff also observed that maximal Hausdorff fields have uncountable cofinality; indeed, he proved the stronger result that the underlying ordered set of a maximal Hausdorff field $H$ is   $\eta_1$: {\it if
$A, B$ are countable subsets of $H$ and $A<B$, then $A<h<B$ for some~${h\in H}$.}\/
A real closed ordered field is $\aleph_1$-saturated iff its underlying ordered set is~$\eta_1$. Standard facts from model theory (or \cite{EGH}) now yield an observation that could have been made by Hausdorff himself in the wake of Artin and Schreier~\cite{AS}:

\begin{cor}\label{cor:Hausdorff}
Assuming \textup{CH} \textup{(}the Continuum Hypothesis\textup{)}, all maximal Hausdorff fields are isomorphic. 
\end{cor}

\noindent
This observation was in fact made by Ehrlich~\cite{Ehr12} in the more specific form that under CH any maximal Hausdorff field is isomorphic to the field of surreal numbers of countable length; see Section 5 below for basic facts on surreals. 
We don't know whether here the assumption of \textup{CH} can be omitted.
(By~\cite{Esterle}, the negation of \textup{CH} implies the existence of non-isomorphic real closed $\eta_1$-fields of size~$2^{\aleph_0}$.) It may also be worth mentioning that the intersection of all maximal Hausdorff fields is quite small: it is just the field of real algebraic numbers. 

\subsection*{Hardy fields}

\noindent
A Hardy field is a Hausdorff field whose germs can be differentiated. 
This leads to a much richer theory. To define Hardy fields formally we introduce the subring
$$\mathcal C^n := \{ f\in\mathcal C: \text{$f$ is eventually $n$ times continuously differentiable} \}$$
of $\mathcal C$, with $\mathcal C^0=\mathcal C$. Then each $f\in\mathcal C^{n+1}$ has derivative $f'\in\mathcal C^n$.
A {\it Hardy field\/} is a subfield of~$\mathcal{C}^1$ that is closed under $f\mapsto f'$;
Hardy fields are thus not only ordered fields but also differential fields.
The Hausdorff fields listed in \eqref{ex:Hausdorff fields} are all Hardy fields;
moreover, for each o-minimal expansion $\tilde{\R}$ of the field of reals, $H(\tilde{\R})$ is a Hardy field.
As with Hausdorff fields,  each Hardy field is contained in a maximal one.
For an element $f$ of a Hardy field we have either $f'>0$, or $f'=0$, or $f'<0$, so $f$ is
either eventually strictly increasing, or eventually constant, or eventually strictly decreasing.
(This may fail for $f$ in a Hausdorff field.)
Each element of a Hardy field  is contained in the intersection $\bigcap_n \mathcal C^n$, but
not necessarily in its  subring~$\mathcal C^\infty$ consisting of those germs which are eventually infinitely differentiable. 
%Some authors (like \cite{Sj71}) focus on Hardy fields which are contained in $\mathcal C^\infty$, or in its subring
%$$\mathcal C^\omega := \{ f\in\mathcal C: \text{$f$ is eventually analytic} \}.$$
In a Hardy field $H$, the ordering and derivation interact in a pleasant way: if $f\in H$ and $f>\R$, then $f'>0$. Asymptotic relations in $H$ can be differentiated and integrated:
for $0\neq f,g\nasymp 1$ in $H$,  $f\preceq g$ iff $f'\preceq g'$.

\subsection*{Extending Hardy fields}

Early work on Hardy fields focussed on solving algebraic equations
and simple first order differential equations: Borel~{\cite{Bor99}}, Hardy~{\cite{Har12a,Har12b}}, Bourbaki~{\cite{Bour51}}, Mari\'c~\cite{Mar72},
Sj\"odin~{\cite{Sj71}}, Robinson~\cite{Rob}, Rosenlicht~\cite{Ros83}.  As a consequence, every Hardy field~$H$
has a smallest real closed Hardy field extension $\operatorname{Li}(H)\supseteq \R$ that is also closed under integration and
exponentiation; call $\operatorname{Li}(H)$ the {\em Hardy-Liouville closure\/} of~$H$.  
(Hardy's field of   LE-functions mentioned earlier is contained in~$\operatorname{Li}(\R)$.)  Here is a rather general result of this kind, due to Singer~{\cite{Sing75}}:

\begin{theorem}
  \label{Singer-th} If $y\in \Cc^1$ satisfies a differential equation $y'P(y)=Q(y)$ where
  $P(Y)$ and $Q(Y)$ are polynomials over a Hardy field $H$
and $P(y)$ is a unit in $\Cc$, 
then~$y$ generates a Hardy field $H \langle y \rangle=H(y,y')$ over $H$.
\end{theorem}

\noindent
Singer's theorem clearly does not extend to second order differential
equations: the nonzero solutions of $y'' + y = 0$ in $\Cc^2$ do not belong to any Hardy
field. The solutions in 
$\Cc^2$ of the  differential equation 
\begin{eqnarray}
    y'' + y & = & \ex^{x^2}  \label{Bos-eq}
\end{eqnarray}
form a two-dimensional affine space $y_0+\R\sin x + \R\cos x$
over $\R$, with $y_0$ any particular solution. Boshernitzan~{\cite{Bos87}} proved that any
of these continuum many solutions generates a Hardy field. Since no Hardy field can contain more than one
solution, there are at least continuum many different maximal Hardy fields.  By the above, each of them contains $\R$, is
real closed, and closed under integration and exponentiation. 
{\it What more can we say about maximal Hardy fields?}\/ 
To give an answer to this question, consider the following   conjectures about Hardy fields~$H$:
\begin{itemize}
\item[A.] {\it For any differential polynomial $P(Y)\in H\{Y\}=H[Y,Y',Y'',\dots]$ and $f<g$ in $H$ with $P(f)<0<P(g)$ there
exists $y$ in a Hardy field extension of $H$ such that $f<y<g$ and $P(y)=0$.}\/
\item[B.] {\it For any countable subsets $A<B$ in $H$ there exists $y$ in a Hardy field extension of $H$ such that $A<y<B$.}\/
\end{itemize}
Conjecture~A for $P\in H[Y,Y']$ holds by \cite{vdD00}. Conjecture~A implies that {\it all maximal Hardy fields are elementarily equivalent\/} as we shall see in
in Section~\ref{sec:H-fields}. 
%we'll see that this conjecture yields }\/ with an axiomatization of their common theory.
Conjecture~B was first raised as a question by Ehrlich~\cite{Ehr12}.   The conjectures together imply that, under~CH, {\it all maximal Hardy fields are isomorphic}\/ (the analogue of Corollary~\ref{cor:Hausdorff}).
We sketch a program to prove A and B in Section~\ref{HH}.

\subsection*{Transseries}

Hardy made the point that the LE-functions seem to cover all orders of infinity that occur naturally in mathematics~\cite[p.~35]{Har10}. 
%``No function has yet presented itself in analysis the laws of whose increase, in so far as they can be stated at all, cannot be stated, so to say, in logarithmico-exponential terms.'' 
But he also suspected that the order of infinity of the compositional inverse of $(\log x)(\log \log x)$ differs from that of any LE-function~\cite{Har12a}; this suspicion is correct.
For a more revealing view of orders of infinity and a more comprehensive theory we need transseries.
For example, transseries lead to an easy argument to confirm  Hardy's suspicion~\cite{DMM97, vdH:phd}.
 Here we  focus on the field $\T$ of LE-series and 
 in accordance with~\cite{vdH:b},    simply call its elements  {\it transseries,}\/ bearing in mind that   many variants of  formal series, such as those appearing in~\cite{Schm01} (see Section~\ref{HT} below), can also  rightfully be called ``transseries''.

Transseries are formal series $f=\sum_{\mathfrak m} f_{\mathfrak m} \mathfrak m$ 
where the~$f_{\mathfrak m}$ are real coefficients and the $\mathfrak m$ are ``transmonomials'' such as 
$$x^r\ (r\in\R),\quad x^{-\log x}, \quad
\ex^{x^2\ex^x},\quad  \ex^{\ex^{x}}.$$
One can get a sense by considering an example like
$$ 7\ex^{\ex^x+\ex^{x/2}+\ex^{x/4}+\cdots} - 3\ex^{x^2} + 5x^{\sqrt{2}} - (\log x)^{\pi} + 42 +
x^{-1} + x^{-2} + \cdots + \ex^{-x}.$$ 
Here think of $x$ as positive infinite: $x>\R$. The transmonomials in this series are arranged from left to right in decreasing order. The reversed order type of the set of transmonomials that occur in
a given transseries can be any countable ordinal. (In the example above it is $\omega+1$ because of the term $\ex^{-x}$ at the end.)
Formally, $\T$ is an ordered  subfield of a Hahn field~$\R[[G^{\text{LE}}]]$ where~$G^{\text{LE}}$ is  the ordered group  of transmonomials (or LE-monomials).
More generally, let $\fM$ be any (totally) ordered commutative group, multiplicatively written, the $\fm\in \fM$ being thought of as monomials, with the ordering
denoted by $\preceq$. The  Hahn field~$\R[[\fM]]$  consists of the formal series $f=\sum_\fm f_\fm \fm$ with real coefficients $f_{\fm}$ whose support $\supp f:=\{\fm\in \fM:\ f_{\fm}\neq 0\}$ is {\it well-based}\/, that is, well-ordered in the reversed ordering~$\succeq$ of~$\fM$. 
Addition and multiplication of these Hahn series works just as for ordinary power series, and the ordering of
$\R[[\fM]]$ is determined by declaring a nonzero Hahn series to be  positive if its leading coefficient is positive (so the series above, with leading coefficient $7$, is positive). Both~$\R[[G^{\text{LE}}]]$ and its ordered subfield $\T$ are real closed. Informally, each transseries is obtained, starting with the powers $x^r$~($r\in\R$), by applying the following operations finitely many times:
\begin{enumerate}
\item multiplication with real numbers;
\item infinite summation in $\R[[G^{\text{LE}}]]$;
\item exponentiation and taking logarithms of positive transseries.
\end{enumerate}
To elaborate on (2), a family $(f_i)_{i \in I}$ in $\R[[\mathfrak{M}]]$ is said to be {{\it summable}} if for each~$\mathfrak{m}$ there
are only finitely many $i\in I$ with $\mathfrak{m} \in \operatorname{supp} f_i$, and
$\bigcup_{i \in I}
\operatorname{supp} f_i$ is well-based; in this case we define the {{\it sum}}
$f = \sum_{i \in I} f_i\in \R [[\mathfrak{M}]]$ of this family by $f_{\mathfrak{m}} =
\sum_{i\in I} (f_i)_{\mathfrak{m}}$ for
each~$\mathfrak{m}$. One can develop a   ``strong'' linear algebra for this notion of ``strong'' (infinite) summation~\cite{vdH:ln,Schm01}. As to (3), it may be instructive to see how to exponentiate a transseries $f$: decompose~$f$ as $f=g+c+\varepsilon$ where $g:=\sum_{ \fm\succ 1} f_{\fm}\fm$ is the infinite part of
$f$,
$c:=f_1$ is its constant term, and $\varepsilon$ its infinitesimal part (in our example   $c=42$ and $\varepsilon=x^{-1} + x^{-2} + \cdots + \ex^{-x}$); then
$$\ex^f\  =\  \ex^g\cdot \ex^c\cdot \sum_{n}  \frac{\varepsilon^n}{n!}$$
where $\ex^g\in\fM$ is a transmonomial,  and $\ex^c\in\R$, $\sum_{n} \frac{\varepsilon^n}{n!}\in 
\R[[G^{\text{LE}}]]$ have their usual meaning. The story with logarithms is a bit different: taking logarithms may also create transmonomials, such as $\log x$, $\log\log x$, etc.

The formal definition of $\T$ is inductive and somewhat lengthy; see~\cite{DMM01,Edgar,vdH:ln} or \cite[Appendix~A]{vdH:b} for   detailed expositions. 
We only note here that by  virtue of the construction of $\T$, series like $\frac{1}{x}+\frac{1}{\ex^x}+\frac{1}{\ex^{\ex^x}}+\cdots$ or $\frac{1}{x}+\frac{1}{x\log x}+\frac{1}{x\log x\log\log x}+\cdots$ 
(involving ``nested'' exponentials or logarithms of unbounded depth), though they are legitimate elements of
$\R[[G^{\text{LE}}]]$,
do not appear in~$\T$; moreover, the sequence $x,\ex^x,\ex^{\ex^x},\dots$ is cofinal in $\T$, and
the sequence $x,\log x,\log\log x,\dots$ is coinitial in the set~${\{f\in\T:f>\R\}}$.
The map $f\mapsto\ex^f$ is an isomorphism of the ordered additive group of~$\T$ onto its multiplicative group  of positive elements, with inverse~${g\mapsto\log g}$. 
As an ordered exponential field,~$\T$ turns out to be 
an elementary extension of~$\R_{\exp}$~\cite{DMM97}. 

Transseries  can be differentiated termwise; for instance, $\left( \sum_n n!\frac{\ex^x}{x^{n+1}}\right)'=\frac{\ex^x}{x}$. We obtain a derivation $f\mapsto f'$ on the field $\T$
with constant field $\{f\in\T:{f'=0}\}=\R$ and satisfying $(\exp f)'=f'\exp f$ and $(\log g)'=g'/g$ for $f,g\in\T$, $g> 0$. Moreover,
each $f\in\T$ has an antiderivative in $\T$, that is, $f=g'$ for some~${g\in\T}$. As in Hardy fields,  $f>\R\Rightarrow f'>0$, for transseries~$f$.
We also have a dominance relation on $\T$: for $f,g\in\T$ we set
\begin{align*}
f\preceq g &\ :\Longleftrightarrow\ 
 \abs{f}\leq c\abs{g}  \text{ for some $c\in\R^{>0}$ } \\
 &\ \ \Longleftrightarrow\ \text{(leading transmonomial of $f$)}\ \preceq\ \text{(leading transmonomial of $g$),}
\end{align*}
and as in Hardy fields we declare 
$f\asymp g :\Longleftrightarrow f\preceq g$ and $g\preceq f$, as well as
${f\prec g} :\Longleftrightarrow f\preceq g$ and $g\not\preceq f$.
As in Hardy fields we can also differentiate
and integrate asymptotic relations: for $0\neq f,g\nasymp 1$ in $\T$ we have $f\preceq g$ iff $f'\preceq g'$.

Hardy's ordered exponential field of (germs of) logarithmic-exponential functions embeds uniquely into $\T$ so as to preserve real constants and to send the germ~$x$ to the transseries $x$; this embedding also preserves the derivation.
However,
the  field of LE-series enjoys
many   closure properties that the  field of LE-functions lacks. For instance, $\T$  is not only closed under
exponentiation and integration, but also comes with a natural operation of composition:
for $f,g\in\T$ with $g>\R$ we can substitute $g$ for $x$ in $f=f(x)$ to obtain $f\circ g=f(g(x))$. The Chain Rule holds: $(f\circ g)'=(f'\circ g)\cdot g'$.  Every $g>\R$ has a compositional inverse in $\T$: a transseries~${f>\R}$ with $f\circ g=g\circ f=x$.
As shown in~\cite{vdH:ln},  a Newton diagram method can be used to solve any ``feasible'' algebraic differential equation in $\T$ (where
the meaning of feasible can be made explicit).

Thus  it is not surprising that
soon after the introduction of~$\mathbb T$ the idea emerged that it should play the role of a {\it universal domain}\/ (akin to Weil's use of this term in algebraic geometry) for asymptotic differential algebra: that it {\it is 
truly the algebra-from-which-one-can-never-exit and that it marks an almost  
impassable horizon for ``ordered analysis''}\/~\cite[p.~148]{Ec92}.
%One possible way to clarify the relationship between Hardy fields and transseries is via the model theory of $H$-fields. 
Model theory provides a  language to make such an intuition precise,
as we explain in our  survey~\cite{ADHOleron}  where we sketched a program to establish the basic model-theoretic properties of $\T$, carried out in~\cite{vdH:b}.  
Next we briefly discuss our main results from~\cite{vdH:b}.

\section{$\mathbf{H}$-Fields}\label{sec:H-fields}

We shall consider $\T$ as an $\mathcal L$-structure where the language $\mathcal L$ has the primitives
$0$, $1$, $+$, $-$, $\,\cdot\,$, $\der$ (derivation), $\leq$ (ordering), $\preceq$ (dominance).
More generally, let~$K$ be any ordered differential field with constant field $C=\{f\in K:f'=0\}$. This
yields a dominance relation $\preceq$ on $K$ by
$$f\preceq g  \quad :\Longleftrightarrow\quad \text{$\abs{f}\leq c\abs{g}$ for some positive $c\in C$} $$
and we view $K$ accordingly as an $\mathcal L$-structure. The convex hull of $C$ in $K$ is the valuation ring
$\mathcal O = \{ f\in K: f\preceq 1\}$  of $K$, with its maximal ideal $\smallo := {\{ f\in K : f\prec 1\}}$ of infinitesimals. 

\begin{definition}
An {\it $H$-field}\/  is an ordered differential field $K$ such that (with the notations above),
$\mathcal O=C+\smallo$, and for all $f\in K$ we have: $f>C \Rightarrow f'>0$. 
\end{definition}

\noindent
Examples include all Hardy fields that contain $\R$, and all ordered differential subfields of $\T$ that contain $\R$. In particular, $\T$ is an $H$-field, but $\T$ has further basic elementary properties that do not follow from this: its derivation is small, and it is Liouville closed.
An $H$-field $K$ is said to have {\em small derivation\/} if it satisfies $f\prec 1\Rightarrow f'\prec 1$,
and to be {\em Liouville closed\/} if it is real closed and for every $f\in K$ there are
$g,h\in K$, $h\neq 0$, such that $g' = f$ and $h'=hf$. 
Each Hardy field $H$ has small derivation, and $\operatorname{Li}(H)$ is Liouville closed.

Inspired by the familiar characterization of real closed ordered fields via the intermediate value property for
one-variable polynomial functions,
we say that an $H$-field $K$ has the {\it Intermediate Value Property}\/ (IVP) if for all differential polynomials $P(Y)\in K\{Y\}$ and all $f<g$ in $K$ with $P(f)<0<P(g)$ there
is some $y\in K$  with $f<y<g$ and $P(y)=0$.  Van der Hoeven showed that a certain variant of $\T$, namely its $H$-subfield of gridbased transseries, has IVP; see~\cite{vdH:int}. 

\begin{theorem}\label{thm:complete axioms}
The $\mathcal L$-theory of $\T$ is completely axiomatized by the requirements:
being an $H$-field with small derivation;
being Liouville closed; and
having IVP.
\end{theorem}

\noindent
Actually, IVP is a bit of an afterthought:
in~\cite{vdH:b} we use other (but 
equivalent) axioms that will be detailed below.
We mention the above variant for expository reasons and
since it explains why Conjecture~A from Section~\ref{sec:trans} yields
that all maximal Hardy fields are elementarily equivalent. Let us define an {\it $H$-closed field\/} to be
an $H$-field that is Liouville closed and has the IVP.
All $H$-fields embed into $H$-closed fields, and the latter are exactly the
existentially closed $H$-fields. Thus:

\begin{theorem}\label{thm:model complete}
The theory of $H$-closed fields is model complete.
\end{theorem}

\noindent
Here is an unexpected byproduct of our proof of this theorem:

\begin{cor}\label{cor:model complete}  $H$-closed fields have no proper
differentially algebraic $H$-field extensions with the same constant field. 
\end{cor}

\noindent
IVP refers to the ordering, but the valuation given by $\preceq$ is more robust and more useful.
IVP comes from two more fundamental properties: {\it $\upo$-freeness}\/ and {\it newtonianity}\/
(a differential version of  henselianity). These concepts make sense for any
 differential field with a suitable dominance relation $\preceq$ in which the equivalence
$f\preceq g\Longleftrightarrow f'\preceq g'$ holds for   $0\neq f,g\prec 1$.

To give an inkling of these somewhat technical notions, let $K$ be an $H$-field and assume that for every $\phi\in K^\times$ for which the derivation $\phi\der$  is small
(that is, $\phi\der\smallo\subseteq \smallo$), there exists $\phi_1\prec \phi$ in $K^\times$ such that $\phi_1\der$ is small. (This assumption is satisfied for Liouville closed
$H$-fields.) Let $P(Y)\in K\{Y\}^{\neq}$. We wish to understand how the function $y\mapsto P(y)$
behaves for $y\preceq 1$.
It turns out that this function only reveals its true colors after rewriting $P$ in terms of a 
derivation~$\phi\der$ with suitable~$\phi\in K^\times$.

Indeed, this rewritten $P$ has the form 
$a\cdot(N + R)$
with~$a\in K^\times$ and where $N(Y)\in C\{Y\}^{\neq}$ is  independent of~$\phi$ for sufficiently small $\phi\in K^\times$ with respect to~$\preceq$, subject to $\phi\der$ being small, and where the
coefficients of~$R(Y)$ are infinitesimal. We call $N$ the {\it Newton polynomial}\/ of~$P$.
Now~$K$ is said to be {\it $\upo$-free}\/ if for all $P$ as above its Newton polynomial
has the form $A(Y)\cdot(Y')^n$ for some $A\in C[Y]$ and some~$n$. We say that~$K$ is {\it newtonian}\/ if for all $P$ as above with~$N(P)$ of degree~$1$ we have
$P(y)=0$ for some~${y\in\mathcal O}$.  
For $H$-fields,  IVP $\Longrightarrow$ $\upo$-free and newtonian; for Liouville closed $H$-fields, the converse also holds.

Our main result in \cite{vdH:b} refines Theorem~\ref{thm:model complete}  by
giving quantifier elimination for the theory of $H$-closed fields  in the language~$\mathcal L$ above augmented by an additional unary function symbol $\iota$ and two extra
unary predicates $\Upl$ and $\Upo$. These have defining axioms in terms of the other primitives. 
Their
interpretations in $\T$ are as follows:
$\iota(f)=1/f$ if $f\neq 0$, $\iota(0)=0$, and with $\ell_0:=x$, $\ell_{n+1}:=\log\ell_n$,
\begin{align*}
\Upl(f)&\quad\Longleftrightarrow\quad f<\upl_n:=\textstyle\frac{1}{\ell_0}+\frac{1}{\ell_0\ell_1}+\cdots+\frac{1}{\ell_0\cdots\ell_n}\ \text{for some $n$,}\\
\Upo(f)&\quad\Longleftrightarrow\quad f<\upo_n:=\textstyle\frac{1}{(\ell_0)^2}+\frac{1}{(\ell_0\ell_1)^2}+\cdots+\frac{1}{(\ell_0\cdots\ell_n)^2}\ \text{for some $n$.}
\end{align*} 
Thus $\Upl$ and $\Upo$ define downward closed subsets of $\T$. The sequence $(\upo_n)$ also appears in
classical non-oscillation theorems for second-order linear differential equations.
The $\upo$-freeness of $\T$ reflects the fact that  $(\upo_n)$ has no pseudolimit in the valued field~$\T$.
Here are some applications of this quantifier elimination:

\begin{cor}\label{cor:QE} \mbox{}
\begin{enumerate}
\item[\textup{(1)}] ``O-minimality at infinity'': if $S\subseteq\T$ is definable, then
for some $f\in\T$ we either have $g\in S$ for all $g>f$ in $\T$ or
$g\notin S$ for all $g>f$ in $\T$.
%$(f,+\infty)\subseteq X$ or $(f,+\infty)\cap X=\emptyset$.
\item[\textup{(2)}] All  subsets of $\R^n$ definable in $\T$ are semialgebraic.
\end{enumerate}
\end{cor}

\noindent
Corollaries~\ref{cor:model complete} and \ref{cor:QE} are the departure point  for 
developing a notion of (dif\-fer\-en\-tial-algebraic) dimension for definable sets in $\T$; see~\cite{ADHdim}.

The results reported on above make us confident that the category of $H$-fields is the right setting
for  asymptotic differential algebra. To solidify this impression we return to the motivating
examples---Hardy
fields, ordered differential fields of transseries, and surreal numbers---and
consider how they are related.
We start with Hardy fields, which  historically came  first.

\section{$\mathbf{H}$-Field Elements as Germs}\label{HH}

After Theorem~\ref{Singer-th} and Boshernitzan~\cite{Bos82,Bos87},
the first substantial ``Hardy field'' result on more general differential equations was obtained by van der Hoeven~{\cite{vdH:hfsol}}. In what follows we use ``$\d$-algebraic'' to mean ``differentially algebraic'' and  ``$\d$-transcendental'' to mean ``differentially transcendental''.

\begin{theorem}
  \label{hfsol-th} The differential subfield $\T^{\da}$ of $\T$ whose elements are the
   $\d$-al\-ge\-braic
  transseries is isomorphic over $\R$ to a Hardy field.
\end{theorem}

\noindent
The proof of this theorem is in the spirit of model theory, 
iteratively extending  by a single $\d$-algebraic transseries.
The most difficult case (immediate extensions) is handled through
  careful construction of suitable solutions as convergent series of
iterated integrals.
We are currently trying to generalize
Theo\-rem~\ref{hfsol-th} to $\d$-algebraic extensions of arbitrary Hardy fields.
Here is our plan:
%Our plan is as follows:

\begin{theorem}\label{upofree}
Every Hardy field has
an $\upomega$-free Hardy field extension.
\end{theorem}

\begin{theorem}[in progress]\label{dalg-gen-th}
Every $\upomega$-free Hardy field has a newtonian
$\d$-al\-ge\-braic Hardy field extension.
\end{theorem}

\noindent
These two theorems, when established, imply that all maximal Hardy fields are 
$H$-closed. Hence (by Theorem~\ref{thm:complete axioms}) they will 
all be elementarily equivalent to $\T$, and since $H$-closed fields have the IVP,  
Conjecture~A from Section~\ref{sec:trans} will follow.
 
 In order to get an even better grasp on the structure of maximal Hardy fields,
we also need to understand how to adjoin $\d$-transcendental germs to Hardy fields. An example of this situation is
given by $\d$-transcendental series such as
$\sum_{n} n!!x^{- n}$. By an old result by \'E.~Borel~\cite{Bor95} every
formal power series $\sum_n a_nt^n$ over~$\R$ is the Taylor series at $0$ of a
$\Cc^{\infty}$-function $f$ on $\R$; then $\sum_n a_nx^{-n}$ is an
asymptotic expansion of
the function $f(x^{-1})$ at $+\infty$, and it is easy to show that if this series is $\d$-transcendental, then
the germ at $+\infty$ of this function does generate a Hardy field.  Here is a far-reaching generalization:

\begin{theorem}[in progress]
\label{imm-Hardy-th} 
Every pseudocauchy sequence $(y_n)$ in a Hardy field~$H$ has a pseudolimit in some Hardy field extension of
$H$.
\end{theorem}

\noindent
The proof of this for $H$-closed $H\supseteq \R$ relies heavily on results from \cite{vdH:b}, using also intricate 
 glueing techniques. For extensions that increase the value group, we need very different constructions. If succesful, 
 these constructions  in combination with Theorem~\ref{imm-Hardy-th} will lead to a proof of 
 Conjecture~B from Section~\ref{sec:trans}:

\begin{theorem}[in progress]
\label{intermediate-Hardy-th} For any countable subsets $A < B$ of a
Hardy field~$H$ there exists an element $y$ in a Hardy field extension of $H$
with $A < y < B$.
\end{theorem}

\noindent
The case $H\subseteq\mathcal C^\infty$, $B=\emptyset$ was already dealt with by Sj{\"o}din~{\cite{Sj71}}. 
%The general statement was first raised as a question by Ehrlich~{\cite{Ehr12}}.
The various ``theorems in progress'' together with results from
\cite{vdH:b} imply that any maximal Hardy fields $H_1$ and $H_2$
are back-and-forth equivalent, which 
is considerably stronger than~$H_1$ and~$H_2$  being elementarily
equivalent. It implies for example
\begin{quote}
\textit{Under ${\rm CH}$ all maximal Hardy fields are isomorphic}.
\end{quote}
This would be the Hardy field analogue of Corollary~\ref{cor:Hausdorff}. (In contrast to maximal Hausdorff fields, however, maximal Hardy fields cannot be
$\aleph_1$-saturated, since their constant field is $\R$.) 
When we submitted this manuscript, we had finished the proof of Theorem~\ref{upofree}, and also the proof of Theorem~\ref{imm-Hardy-th} in the relevant $H$-closed case. 

\subsection*{Related problems} Some authors (such as \cite{Sj71}) prefer to consider only
Hardy fields contained in~$\Cc^{\infty}$. Theorem~\ref{upofree}
and our partial result for Theorem~\ref{imm-Hardy-th} go through 
in the~$\Cc^{\infty}$-setting. 
All the above ``theorems in progress'' are plausible in that setting.

What about real analytic Hardy fields (Hardy fields contained in the subring~$\Cc^{\omega}$ of $\mathcal C$
consisting of all real analytic germs)?
In that setting Theorem~\ref{upofree} goes through. Any $\d$-algebraic Hardy field extension
of a real analytic Hardy field is itself real analytic, and so Theorem~\ref{dalg-gen-th} (in progress)
will hold in that setting as well. 
However, our glueing technique employed in
the proof of Theorem~\ref{imm-Hardy-th} doesn't work there. 
 
Kneser~{\cite{Kne50}} obtained a real
analytic solution $E$ at infinity to the functional equation $E (x + 1) = \exp E (x)$. It grows faster than any finite iteration of the exponential function,
and   generates a Hardy field. See Boshernitzan~{\cite{Bos86}} for results of this kind, and a
proof that  Theorem~\ref{intermediate-Hardy-th} holds for $B=\emptyset$ in the real analytic setting. 
So  in this context we also have an abundant supply of Hardy fields.

Similar issues arise for germs of quasi-analytic and 
``cohesive'' functions~\cite{Ec92}. These classes of functions are somewhat more flexible
than the class of real analytic functions. For instance, the series $x^{- 1} + \ex^{- x} + \ex^{- \ex^x} + \cdots$ converges uniformly
for~${x>1}$ to a~cohesive function that is not real analytic.

\subsection*{Accelero-summation}

The definition of a Hardy field ensures that the differential field
operations never introduce oscillatory behavior. Does this behavior persist for operations such as
composition or various integral transforms? In this connection we note that
the  Hardy field $H(\tilde{\R})$ associated to an o-minimal expansion $\tilde{\R}$ of the field of reals
 is always closed under composition (see~\cite{Miller}). 

To illustrate the problem with composition, let $\alpha$ be a real number $>1$ and let $y_0\in \Cc^2$ be a solution to~\eqref{Bos-eq}. Then $z_0:=y_0 (\alpha x)$ satisfies the  
equation
\begin{equation}
  \alpha^{- 2} z'' + z \ = \ \ex^{\alpha^2 x^2}.  \label{Bos-eq2}
\end{equation}
It can be shown that $\{y_0 + \sin x,  z_0\}$
generates a Hardy field, but it is clear that no Hardy
field containing both $y_0 + \sin x$ and $z_0$ can be closed under composition.

Adjoining  solutions
to~\eqref{Bos-eq} and~\eqref{Bos-eq2} ``one by one'' as in the proof of Theorem~\ref{hfsol-th} will not
prevent the resulting Hardy fields to contain both $y_0 + \sin x$ and~$z_0$. In order to obtain closure under composition
we therefore need an alternative device. {\'E}calle's theory of
{\it accelero-summation}~{\cite{Ec92}} is much more than  that. Vastly extending
Borel's  summation method for divergent series~\cite{Bor99}, it associates to each {\em accelero-summable\/} transseries an
{\em analyzable\/} function. In this way many non-oscillating
real-valued functions that arise naturally (e.g., as solutions of algebraic differential equations) can be represented faithfully by transseries. 
This leads us to
conjecture an improvement on
Theorem~\ref{hfsol-th}:

\begin{conjecture}
  \label{accsum-conj}   Consider the real accelero-summation process
  where we systematically use the organic average whenever we encounter
  singularities on the positive real axis. This
  yields a composition-preserving $H$-field isomorphism  from $\T^{\da}$
  onto a Hardy field contained in $\mathcal C^{\omega}$.  
\end{conjecture}

\noindent
There is little doubt that this  holds. The main difficulty 
here is that a~full proof will involve many tools forged
by {\'E}calle in connection with
accelero-summation, such as resurgent functions, well-behaved averages,
cohesive functions, etc., with some of these tools
requiring further elaboration; see also~{\cite{Cos08,Men96}}.

The current theory of
accelero-summation only sums transseries with coefficients
in $\R$. Thus it is not clear how to generalize
Conjecture~\ref{accsum-conj} in the direction of Theorem~\ref{dalg-gen-th}.
Such a generalization might require introducing transseries over a Hardy
field~$H$ with suitable additional structure, as well as a corresponding
theory of accelero-summation over $H$ for such transseries. In particular, elements of
$H$ should be accelero-summable over $H$ in this theory, by construction.

\section{$\mathbf{H}$-Field Elements as Generalized Transseries}\label{HT}

\noindent
Next we discuss when $H$-fields embed into differential fields of formal series.  
A classical embedding theorem of this type is due to Krull~\cite{Krull32}:  any
valued field has a spherically complete immediate extension.
As a consequence, any real closed field containing $\R$ is isomorphic over $\R$ to a subfield of a Hahn field~$\R[[\mathfrak{M}]]$
with divisible monomial group
$\mathfrak{M}$, such that the subfield contains~$\R(\mathfrak{M})$.
We recently proved an analogue of this   theorem for valued differential fields~{\cite{vdH:vdf}}. Here a 
 {{\it valued differential field}} is a valued field   of
equicharacteristic zero equipped with a derivation   that is continuous with respect to the
valuation topology. 
%We denote the corresponding valuation ring by $\mathcal{O}=\mathcal{O}_K$ and the corresponding maximal ideal of $\mathcal{O}$ by $\smallo=\smallo_K$. The continuity requirement gives the existence of an element $\phi \in K$ with $\der \smallo \subseteq \phi \smallo$. A valued differential field extension $L \supseteq K$ is said to be {{\it strict}} if for every $\phi \in K$, 
%$$  \der \smallo_K \subseteq \phi \smallo_K \ \Longrightarrow \  \der \smallo_L \subseteq \phi \smallo_L,\qquad  \der \mathcal{O}_K \subseteq \phi \mathcal{O}_K \ \Longrightarrow \   \der \mathcal{O}_L \subseteq \phi \mathcal{O}_L.$$

\begin{theorem}
  \label{diff-Krull-th}Every valued differential field has a
spherically complete immediate   extension.
\end{theorem}

\noindent
For a real closed $H$-field $K$ with constant field $C$ this theorem gives a Hahn field $\widehat{K}=C[[\mathfrak{M}]]$ with a derivation ${\der}$
on $\widehat{K}$ making it an $H$-field with constant
field~$C$ such that~$K$ is isomorphic over $C$ to an $H$-subfield
of $\hat{K}$ that contains~$C(\mathfrak{M})$.  
A shortcoming of this result is that there is no guarantee that ${\der}$ preserves infinite summation. In contrast, the derivation of
$\T$ is {\it strong} (does preserve
infinite summation). 
An abstract framework for even
more general notions of transseries is due to van der Hoeven and his former
student Schmeling~{\cite{Schm01}}.

\subsection*{Fields of transseries}\label{transseries-field-sec}

To explain this,
consider an (ordered) Hahn field $\R
[[\mathfrak{M}]]$   with a~partially defined function $\exp$
obeying the usual rules of exponentiation;
 see~{\cite[Section~4.1]{vdH:ln}} for details. In particular, $\exp$ has  a partially defined inverse function~$\log$. We say that~$\R
[[\mathfrak{M}]]$ is a {{\it field of transseries}} if the
following conditions hold:
\begin{itemize}
  \item[(T1)] the domain of the function $\log$ is  $\R [[\mathfrak{M}]]^{>0}$;
  \item[(T2)] for each $\fm\in\fM$ and $\fn\in\supp\log\fm$ we have $\fn\succ 1$;  
  \item[(T3)] $\log (1+\varepsilon) = \varepsilon - \frac{1}{2} \varepsilon^2 +
  \frac{1}{3} \varepsilon^3 + \cdots$ for all $\varepsilon \prec 1$  in 
  $\R[[\mathfrak{M}]]$; and
  \item[(T4)] for every sequence $(\mathfrak{m}_n)$ in $\mathfrak{M}$   with $\mathfrak{m}_{n + 1} \in \operatorname{supp}
  \log \mathfrak{m}_n$ for all $n$, there exists an index $n_0$  such that for all $n\geq n_0$ and all $\mathfrak{n} \in \operatorname{supp}
    \log \mathfrak{m}_n$, we have $\mathfrak{n} \succcurlyeq \mathfrak{m}_{n + 1}$ and
     $(\log \mathfrak{m}_n)_{\mathfrak{m}_{n + 1}} = \pm 1$.
\end{itemize}
The first three axioms record basic facts from the standard construction of
trans\-series. The fourth axiom is more intricate and puts limits on the kind of
``nested transseries'' that are allowed. Nested transseries such as
\begin{eqnarray}
  y & = & \sqrt{x} + \ex^{\sqrt{\log x} + \ex^{\sqrt{\log \log x} +
  \ex^{\cdots}}}  \label{nested-transseries}
\end{eqnarray}
are naturally encountered as solutions of functional equations, in this case
\begin{eqnarray}
  y (x) & = & \sqrt{x} + \ex^{y (\log x)} .  \label{nested-eq}
\end{eqnarray}
Axiom~(T4) does allow nested transseries as in~\eqref{nested-transseries}, but excludes series like
%stipulates that, from a certain point on, we only
%allow terms at the left of infinite diagonals of raising imbricated exponentials. For instance,
\begin{eqnarray*}
  u & = & \sqrt{x} + \ex^{\sqrt{\log x} + \ex^{\sqrt{\log \log x} +
  \ex^{\cdots} + \log \log \log x} + \log \log x} + \log x,
\end{eqnarray*}
%is not considered to be a valid transseries. 
which solves the functional equation $u (x) = \sqrt{x} + \ex^{u (\log x)} +
\log x$; in some sense, $u$ is a  perturbation of the solution $y$ in~\eqref{nested-transseries}
to the equation~\eqref{nested-eq}.
%there is way to express this solution  in terms of~$y$ .

Schmeling's thesis~{\cite{Schm01}} shows how to extend a given field of transseries $K=\R [[\mathfrak{M}]]$ with new
exponentials and nested transseries like~\eqref{nested-transseries}, and if~$K$ also comes with a~strong derivation, 
 how to extend this derivation as well. Again,~(T4) is
  crucial for this task: naive termwise differentiation leads to
a huge infinite sum that turns out to be summable by~(T4). A
{{\it transserial derivation}} is a~strong derivation on~$K$ such that nested
transseries are differentiated in this way. Such a transserial derivation is
uniquely determined by its values on the {{\it log-atomic}} elements: those  
$\lambda\in K$ for which $\lambda, \log \lambda, \log \log \lambda, \ldots$ are all transmonomials in~$\mathfrak{M}$.

We can now state a transserial analogue of Krull's
theorem. This analogue is a consequence of Theorem~\ref{surreal-embed-th}
below, proved in~{\cite{vdH:bm}}.

\begin{theorem}
  \label{trans-embed-th} Every $H$-field with small derivation and
  constant field $\R$ can be embedded over $\R$ 
  into a~field of transseries with  transserial
  derivation.
\end{theorem}

\noindent
For simplicity, we restricted ourselves to transseries over $\R$. The
theory naturally generalizes to transseries over ordered exponential
fields~{\cite{vdH:ln,Schm01}} and it should be possible to extend
Theorem~\ref{trans-embed-th} likewise.

\subsection*{Hyperseries}

Besides derivations, one can also define a notion of composition for
generalized trans\-series~{\cite{vdH:phd,Schm01}}. Whereas certain functional
equations such as~\eqref{nested-eq} can still be solved using nested
transseries, solving the equation
$ E (x + 1)  =  \exp E (x)$
where~$E(x)$ is the unknown, requires extending $\T$ to a field of transseries with composition containing an element $E(x)=\exp_{\omega} x> \T$, called the {{\it iterator}} of~$\exp x$. 
Its compositional inverse 
$\log_{\omega} x$ should then satisfy
$   \log_{\omega} \log x =  (\log_{\omega} x) - 1$,
%\begin{eqnarray*}
%   \log_{\omega} \log x & = & (\log_{\omega} x) - 1,
%\end{eqnarray*}
providing us with a primitive for  $(x \log x \log_2 x
\cdots)^{- 1}$:
\begin{eqnarray*}
  \log_{\omega} x & = & \int \frac{\text{d}x}{x \log x \log_2 x \cdots} .
\end{eqnarray*}
It is convenient to start with iterated logarithms rather than iterated exponentials, and
to introduce  transfinite
iterators~$\log_{\alpha} x$  recursively using
\begin{eqnarray*}
  \log_{\alpha} x & = & \int \frac{\text{d}x}{\prod_{\beta<\alpha}\log_{\beta} x}  \qquad (\alpha \text{ any ordinal}).
\end{eqnarray*}
By {\'E}calle~{\cite{Ec92}} the iterators $\log_{\alpha} x$
with $\alpha < \omega^{\omega}$ and their compositional in\-ver\-ses~$\exp_{\alpha} x$ suffice to
resolve all pure composition equations of the form
$$  f^{\circ k_1} \circ \phi_1 \circ \cdots \circ f^{\circ k_n} \circ \phi_n \ =
  \ x \quad\text{where $\phi_1, \ldots, \phi_n\in\T$ and $k_1, \ldots, k_n
\in \N$.}$$
The resolution of more complicated functional equations
involving differentiation and composition requires the introduction of fields
of {{\it hyperseries}}: besides exponentials and logarithms, hyperseries are
allowed to contain iterators~$\exp_{\alpha} x$ and~$\log_{\alpha} x$ of any
strength~$\alpha$. For~{$\alpha < \omega^{\omega}$}, the necessary
constructions were carried out in~{\cite{Schm01}}. The
ultimate objective is to construct a field $\textbf{Hy}$ of hyperseries
as a proper class, similar to the field of surreal numbers, endow it
with its canonical derivation and composition, and establish the following:

\begin{conjecture}
  Let $\Phi$ be any partial function from $\textbf{Hy}$ into itself,
  constructed from elements in $\textbf{Hy}$,
  using the field operations, differentiation and composition. Let $f < g$
  be   hyperseries in~$\textbf{Hy}$ such that $\Phi$ is defined on the
  closed interval $[f, g]$ and   $\Phi (f) \Phi (g) < 0$. Then for some $y\in\textbf{Hy}$ we have
  $\Phi (y) = 0$ and $f < y < g$.
\end{conjecture}

\noindent
One might then also consider $H$-fields with an additional
composition operator and try to prove that these structures can always be embedded into
$\textbf{Hy}$.

\section{Growth Rates as Numbers}\label{HN}

\noindent 
Turning to surreal numbers, how do they fit into asymptotic differential algebra?

\subsection*{The $\mathbf{H}$-field of surreal numbers}
%Surreal numbers were invented by Conway~{\cite{Con76}}. 
The totality $\textbf{No}$
 of surreal numbers is not a set but a proper class: a surreal $a
\in \textbf{No}$ is uniquely represented by a transfinite sign sequence
$(a_{\lambda})_{\lambda < \ell (a)} \in \{ -, + \}^{\ell (a)}$ where~$\ell (a)$
is an ordinal, called the {{\it length}} of $a$; a surreal $b$ is said to be {\em simpler than\/}~$a$ (notation: $b<_s a$) if
the sign sequence of $b$ is a proper initial segment of that of~$a$.  
Besides the (partial) ordering $<_s$, $\textbf{No}$ also carries a natural (total) lexicographic ordering $<$. For any {\it sets}\/ $L< R$ of surreals there is a unique simplest surreal $a$ with $L < a < R$; this $a$ is denoted
by $\{ L \divides R \}$ and called the {{\it simplest}} or {{\it earliest}}
surreal between $L$ and $R$. In particular, $a = \{L_a \divides R_a \}$ for
any $a \in \textbf{No}$, where $L_a: = \{ b<_s a:\  b < a\}$ and
$R_a = \{ b<_s a:\  b>a\}$. We let~$a^L$ range over elements of $L_a$, and $a^R$ over elements of $R_a$.  

A rather magical property of surreal numbers is that various
operations have natural inductive definitions. For instance, we have ring operations given by
\begin{align*}
  a + b &\ := \ \{ a^L + b, a + b^L \divides a^R + b, a + b^R \}\\
  ab &\ :=\  \{ a^L b + ab^L - a^L b^L, a^R b + ab^R - a^R b^R \divides
  \nosymbol \nobracket\\
    & \phantom{\ \ :=\  \{ \nobracket} \nobracket a^L b + ab^R - a^L b^R, a^R b +
  ab^L - a^R b^L \}.
\end{align*}
Remarkably, these operations make $\textbf{No}$ into a real closed field with $<$ as its field ordering and with $\R$ uniquely embedded as an initial subfield. (A set $A\subseteq \textbf{No}$
is said to be {\em initial\/} if for all $a\in A$ all $b<_s a$ are also in $A$.) 

Can we use such magical recursions to introduce
other reasonable operations? Exponentiation  was dealt with by Gonshor~{\cite{Gon86}}. But it remained long open how to define a ``good'' derivation~$\der$ on 
$\textbf{No}$ such that $\der (\omega) = 1$. (An ordinal $\alpha$ is identified with the surreal of length $\alpha$   whose sign sequence has just plus signs.) 
A positive answer was given recently
by Berarducci and Mantova~{\cite{BM15}}.
Their construction goes in two parts. They first analyze 
$\textbf{No}$ as an exponential field, and show that it is basically a field of transseries in the sense of
Section~\ref{transseries-field-sec}. A transserial
derivation on~$\textbf{No}$ is determined by its values at
log-atomic elements. There is some flexibility here, but {\cite{BM15}} presents a ``simplest'' way to choose these
derivatives. Most important, that choice indeed leads to a derivation~$\der_{\operatorname{BM}}$ on~$\textbf{No}$. In addition:

\begin{theorem}[Berarducci-Mantova~{\cite{BM15}}] \label{BM-th}
The derivation~$\der_{\operatorname{BM}}$ is transserial and makes
$\textbf{\textup{No}}$ a Liouville   closed $H$-field with constant field~$\R$.
%\marginpar{shortened}
%The field  $\textbf{\textup{No}}$ of surreals is a field of transseries. The derivation~$\der_{\operatorname{BM}}$ is transserial   and makes $\textbf{\textup{No}}$ a Liouville   closed $H$-field with constant field~$\R$.
\end{theorem}

\noindent
This result was further strengthened in~{\cite{vdH:bm}}, using key results from~{\cite{vdH:b}}:

\begin{theorem}
  \label{BM-gen-th}  $\textbf{\textup{No}}$  with the derivation $\der_{\operatorname{BM}}$ is an $H$-closed field. %In particular, it has the \textup{IVP}.
\end{theorem}

\subsection*{Embedding $\mathbf{H}$-fields into $\textbf{No}$}

In the remainder
of this section we consider $\textbf{No}$ as equipped with the derivation~$\der_{\operatorname{BM}}$, although Theorems~\ref{BM-th} and~\ref{BM-gen-th} and much of what follows hold for
other transserial derivations. 
Returning to our main topic of embedding $H$-fields into specific $H$-fields such
as~$\textbf{No}$, we also proved the following in~{\cite{vdH:bm}}:

\begin{theorem}
  \label{surreal-embed-th} Every $H$-field with small derivation and
  constant field $\R$ can be embedded as
  an ordered differential field into $\textbf{\textup{No}}$.
\end{theorem}

\noindent
How ``nice'' can we take the embeddings in
Theorem~\ref{surreal-embed-th}? For instance, when can we arrange the image of the embedding to be initial? The image of the natural embedding $\T \to \textbf{No}$ is indeed initial, as has been
shown by Elliot Kaplan. 

For further discussion it is convenient to introduce, given an ordinal $\alpha$, the
{\it set}\/ $\textbf{No} (\alpha):=\big\{a\in\textbf{No}:\ell(a)<\alpha\big\}$.  It turns out that for uncountable cardinals~$\kappa$, $\textbf{No} (\kappa)$ is
closed under the differential field operations, and in \cite{vdH:bm} we also show:

\begin{theorem}
  The $H$-subfield $\textbf{\textup{No}} (\kappa)$ of $\textbf{\textup{No}}$ is an elementary submodel of
  $\textbf{\textup{No}}$.
\end{theorem}

\noindent
In particular, the $H$-field $\textbf{\textup{No}} (\omega_1)$ of surreal numbers of countable length
  is an elementary submodel of $\textbf{\textup{No}}$.
It has the $\eta_1$-property: for any countable subsets~${A < B}$ of~$\textbf{\textup{No}} (\omega_1)$ there exists $y\in \textbf{\textup{No}} (\omega_1)$
with $A < y < B$. This fact and the various ``theorems in progress'' from Section~\ref{HH} imply: 
\begin{quote} \textit{Under ${\rm CH}$ all maximal Hardy fields are isomorphic
  to $\textbf{\textup{No}} (\omega_1)$}.
\end{quote}
This would be an analogue of Ehrlich's observation about maximal Hausdorff fields. 

\subsection*{Hyperseries as numbers and {\it \textbf{vice versa}}}

The similarities in the constructions of the field
of hyperseries $\textbf{Hy}$ and the field of surreal numbers $\textbf{No}$ led van der Hoeven~{\cite[p.~6]{vdH:ln}} to the following:

\begin{conjecture}
  \label{HyNo-conj} There is a natural isomorphism between $\textbf{Hy}$ and $\textbf{No}$ that associates to any hyperseries $f(x) \in \textbf{Hy}$
its value $f (\omega) \in \textbf{No}$.
\end{conjecture}

\noindent
The problem is to make sense of the value of a hyperseries at $\omega$. 
Thanks to Gonshor's exponential function, it is clear how to evaluate ordinary
transseries at~$\omega$. The difficulties start as soon as we wish to
represent surreal numbers that are not of the form $f(\omega)$ with $f(x)$ an ordinary transseries. That is where the iterators~$\exp_{\omega}$ and~$\log_{\omega}$ come into play:
\begin{eqnarray*}
 \exp_{\omega} \omega& := &  \{ \omega, \exp \omega, \exp_2 \omega, \ldots \divides \nosymbol \} 
 \\
 \log_{\omega} \omega& := &  \{ \R \divides \ldots, \log_2 \omega, \log \omega, \omega \} 
 \\
  \exp_{1 / 2} \omega& := &\exp_{\omega} \left( \log_{\omega}  \left( \omega + \tfrac{1}{2} \right)
  \right)\\
  &:= & \left\{ \omega^2, \exp \log^2 \omega, \exp_2 \log_2^2 \omega, \ldots
  \,\middle|\, \ldots, \exp_2  \sqrt{\log \omega}, \exp \sqrt{\omega} \right\}  
\end{eqnarray*}
The intuition behind Conjecture~\ref{HyNo-conj} is that all ``holes in
$\textbf{No}$ can be filled'' using suitable nested hyperseries and suitable iterators
of $\exp$ and $\log$.
It reconciles two {{\it a priori}} very different types of
infinities: on the one hand, we have growth orders corresponding to smooth
functional behavior; on the other side, we have numbers. 
Being able to switch between functions (more precisely: formal series 
acting as functions) and numbers, we may
also transport any available structure in both directions: we immediately
obtain a canonical derivation $\der_{\operatorname{c}}$ (with constant field~$\R$)
and composition $\circ_{\operatorname{c}}$
on $\textbf{No}$, as well
as a notion of simplicity on~$\textbf{Hy}$. 

Does the derivation $\der_{\operatorname{BM}}$  coincide
with the canonical derivation $\der_{\operatorname{c}}$
induced by the conjectured isomorphism? A key observation
is that any derivation~$\der$ on~$\textbf{No}$ with a distinguished
right inverse $\der^{- 1}$ naturally gives rise to a definition of~$\log_{\omega}$:
\begin{align*} \log_{\omega} a \  &:= \ \der^{-1} (\der a \log'_{\omega} a)
\  \text{ where}\\
\log'_{\omega} a &:= 1 \bigg/ \prod_n \log_n a  \qquad (a\in \textbf{No},\ a>\R).
%\log'_{\omega} a &:= 1 / (a \log a \log_2 a \cdots)\qquad (a\in \textbf{No},\ a>\R).
\end{align*} 
(For a family $(a_i)$ of positive surreals, $\prod_i a_i:= \exp \sum_i \log a_i$ if $\sum_i  \log a_i$
is defined.) 
%, thanks to Gonshor, 
%we define $\log'_{\omega} a$ as the value of
%the transseries $\log'_{\omega} x := 1 / (x \log x \log_2 x \cdots)$ at $a$.
Since $\der_{\operatorname{BM}}$  is transserial, it does
admit a distinguished right in\-verse~$\der_{\operatorname{BM}}^{- 1}$. 
According to \cite[Remark~6.8]{BM15}, $\der_{\operatorname{BM}} \lambda = 
1 /\log'_{\omega}\lambda$ for log-atomic $\lambda$ with
$\lambda >\exp_n \omega$ for all~$n$.  For $\lambda=\exp_{\omega}\omega$ and setting $\exp'_{\omega}(a):=\prod_n \log_n \exp_{\omega} a$ for $a\in \textbf{No}^{>0}$,
this yields $\der_{\operatorname{BM}} \lambda = \exp_{\omega}' \omega$, which is also the value we expect for
$\der_{\operatorname{c}}\lambda$.  However, for $\lambda=\exp_{\omega} (\exp_{\omega}\omega)$
we get  $\der_{\operatorname{BM}} \lambda =
\exp_{\omega}'( \exp_{\omega} \omega)$ whereas we expect
$\der_{\operatorname{c}}\lambda= 
(\exp_{\omega}' \omega)\cdot \exp_{\omega}'(\exp_{\omega} \omega)$.
Thus the ``simplest'' derivation $\der_{\operatorname{BM}}$ making $\textbf{No}$ an $H$-field probably does {\it not}\/
coincide with the ultimately~``correct'' derivation
$\der_{\operatorname{c}}$ on $\textbf{No}$.
Berarducci and Mantova~\cite{BM17} use similar considerations to conclude that~$\der_{\operatorname{BM}}$ is incompatible with any reasonable notion of composition
for surreal numbers.

%This was also observed in \cite{BM17} where Berarducci-Mantova formulate desirable properties of a hypothetical notion of composition 
%$(f,g)\mapsto f\circ g\colon\textbf{No}\times\textbf{No}^{>\R}\to\textbf{No}$ 
%for surreal numbers
%and show that~$\der_{\operatorname{BM}}$ is not compatible (in a suitable sense) with any such notion of composition.

\subsection*{The surreal numbers from a model theoretic perspective}

We conclude with  speculations motivated by the  fact that various operations defined by  ``surreal'' recursions 
have a nice model theory.  Examples:
$\left( \textbf{No}; {\leqslant}, {+}, {\,\cdot\,} \right)$ is a model of the
  theory of real closed fields; $\left( \textbf{No}; \leqslant, +, \,\cdot\,, \exp \right)$ is a model of
  the theory of $\R_{\exp}$; and
  $\left( \textbf{No}; {\leqslant}, {+}, {\,\cdot\,}, \der_{\operatorname{BM}} \right)$ is a model
  of the theory of $H$-closed fields.
Each of these theories is model complete in a natural language. Is there a
model theoretic reason that explains why this works so well?

Let us look at this in connection with the last example. Our aim is to define a derivation~$\der$ on
$\textbf{No}$ making it an $H$-field.  Let  $a\in\textbf{No}$ be given for which we wish to define $\der a$,
and assume that $\der b$ has been defined for all $b \in L_a \cup R_a$.
Let~$\Delta_a$ be the class of all surreals $b$ for which there
exists a derivation $\der$ on~$\textbf{No}$ with $\der a= b$ 
and taking the prescribed values on $L_a \cup R_a$. Assembling all conditions that
should be satisfied by $\der a$, it is not hard to see that there exist
sets $L, R\subseteq \textbf{No}$ such that $\Delta_a = \left\{ b \in \textbf{No} : L < b < R
\right\}$. We are left with two main questions:
{\it When do we have $L < R$, thereby
  allowing us to define $\der a \:= \{ L \divides R \}$?}
 {\it Does this lead to a global definition of
  $\der$ on $\textup{\textbf{No}}$ making it an
 $H$-closed field?}\/ 
It might be of interest  to isolate reasonable model theoretic conditions that imply the success
of this type of construction. If the above construction does work, yet
another question is whether the resulting derivation coincides with~$\der_{\operatorname{BM}}$.

\def\bysame{\leavevmode\hbox to3em{\hrulefill}\thinspace}

%setlength{\bibsep}{0pt} % or use whatever dimension you want
%\renewcommand{\bibfont}{\small}

\end{document}